\documentclass[journal,twocolumn,10pt,twoside]{IEEEtran}

\ifCLASSINFOpdf
\else
\fi

\usepackage[utf8x]{inputenc}

\usepackage{hyperref}
\usepackage{amssymb}
\usepackage{graphicx}
\usepackage[cmex10]{amsmath}
\usepackage{array}
\usepackage{mdwtab}
\usepackage{graphics}
\usepackage{epsfig}
\usepackage{times}
\usepackage{url}
\usepackage{lscape}
\usepackage{color}
\usepackage{epsfig}
\newtheorem{theorem}{Theorem}[section]

\newtheorem{remark}[theorem]{Remark}

\usepackage{cite}
\usepackage{balance}
\usepackage{dblfloatfix}
\usepackage{nicefrac}

\definecolor{dark`green'}{rgb}{0,0.50196078431,0}
\definecolor{dark`yellow'}{rgb}{1,0.8,0}

\begin{document}

\title{MPC-Based Emergency Vehicle-Centered Multi-Intersection Traffic Control}

\author{Mehdi Hosseinzadeh, \IEEEmembership{Member, IEEE}, Bruno Sinopoli, \IEEEmembership{Fellow, IEEE}, Ilya Kolmanovsky, \IEEEmembership{Fellow, IEEE}, \\ and Sanjoy Baruah, \IEEEmembership{Fellow, IEEE}  
\thanks{This research has been supported by National Science Foundation under award numbers ECCS-1931738, ECCS-1932530, and ECCS-2020289.}
\thanks{M. Hosseinzadeh and B. Sinopoli are with the Department of Electrical and Systems Engineering, Washington University in St. Louis, St. Louis, MO 63130, USA (email: mehdi.hosseinzadeh@ieee.org; bsinopoli@wustl.edu).}
\thanks{I. Kolmanovsky is with the Department of Aerospace Engineering, University of Michigan, Ann Arbor, MI 48109, USA (email: ilya@umich.edu).}
\thanks{S. Baruah is with the Department of Computer Science and Engineering, Washington University in St. Louis, St. Louis, MO 63130, USA (email:baruah@wustl.edu).}
}


\maketitle

\begin{abstract}
This paper proposes a traffic control scheme to alleviate traffic congestion in a network of interconnected signaled lanes/roads. The proposed scheme is emergency vehicle-centered, meaning that it provides an efficient and timely routing for emergency vehicles. In the proposed scheme, model predictive control is utilized to control inlet traffic flows by means of network gates, as well as configuration of traffic lights across the network. Two schemes are considered in this paper: i) centralized; and ii) decentralized. In the centralized scheme, a central unit controls the entire network. This scheme provides the optimal solution, even though it might not fulfil real-time computation requirements for large networks. In the decentralized scheme, each intersection has its own control unit, which sends local information to an aggregator. The main responsibility of this aggregator is to receive local information from all control units across the network as well as the emergency vehicle, to augment the received information, and to share it with the control units. Since the decision-making in decentralized scheme is local and the aggregator should fulfil the above-mentioned tasks during a traffic cycle which takes a long period of time, the decentralized scheme is suitable for large networks, even though it may provide a sub-optimal solution. Extensive simulation studies are carried out to validate the proposed schemes, and assess their performance. Notably, the obtained results reveal that traveling times of emergency vehicles can be reduced up to $\sim$50\% by using the centralized scheme and up to $\sim$30\% by using the decentralized scheme, without causing congestion in other lanes.
\end{abstract}


\section{Introduction}\label{sec:Introduction}
\IEEEPARstart{T}{raffic} congestion is one of the most critical issues in urbanization. In particular, many cities around the world have experienced 46\textendash70\% increase in traffic congestion \cite{TomTom}. Congested roads not only lead to increased commute times, but also hinder timely deployment of emergency vehicles \cite{Barth2010}. Hence, emergency vehicles often fail to meet their target response time \cite{Oza2021}. According to $\sim$240 million emergency calls every year in the U.S. \cite{911}, such hindering greatly affects hospitalization and mortality rates \cite{Jena2017}.

The common practice by regular vehicles (i.e., non-emergency vehicles) in the presence of an emergency vehicle is to pull over to the right (in two-way roads) or to the nearest shoulder (in one-way roads) \cite{DoT}, and let the emergency vehicle traverse efficiently and timely. This is not always possible, as in dense areas the edges of the roads are usually occupied by parked/moving vehicles.

The chance of an emergency vehicle getting stuck is even higher when it has to traverse intersections with cross-traffic \cite{Hsiao2018}. Note that the majority of incidents involving emergency vehicles happen within intersections \cite{InjuryFacts}. One possible way to cope with this problem is to use traffic lights at intersections to detect emergency vehicles and facilitate their fast and efficient travel. For this purpose, traffic lights in most parts of the U.S. are equipped with proper detectors (e.g., 3M Opticom\textsuperscript{\texttrademark} \cite{Paruchuri2017}), and emergency vehicles are equipped with emitters which broadcast an infrared signal. When the receiver on a traffic light detects a recognized signal, the traffic light changes to allow priority access to the emergency vehicle. In this context, the ``green wave" method has been proposed to reduce emergency vehicles' traveling time \cite{Kang2014}. In the ``green wave" method, a series of traffic lights are successively set to `green' to allow timely passage of emergency vehicles through several intersections \cite{Cao2019}. The main issue with the ``green wave" method is that it leads to prolonged red lights for other lanes \cite{Kapusta2017}, meaning that it may cause congestion in other lanes.

A different method for controlling the traffic in the presence of an emergency vehicle is to convert the traffic control problem to a real-time scheduling problem \cite{Oza2020,Oza2021}. The core idea of this method is to model the vehicles and traffic lights as aperiodic tasks and sporadic servers, respectively, and then to utilize available task scheduling schemes to solve the resulting problem. Other existing traffic control methods either do not consider emergency vehicles \cite{Lin2012,Baskar2012,Tettamanti2014,Jamshidnejad2018,Jafari2019,Rastgoftar2021} or require vehicle to vehicle connectivity \cite{Toy2002,Kamalanathsharma2012,During2014,Weinert2015,Hannoun2019,Wu2020}. Note that the presence of 100\% of connected vehicles is not expected until 2050 \cite{Feng2015}, making these methods inapplicable to the current traffic systems.

The aim of this paper is to propose  control algorithms to manipulate traffic density in a network of interconnected signaled lanes. The core idea is to integrate the Cell Transmission Model (CTM) \cite{Li2017,Shao2018} with Model Predictive Control (MPC) \cite{Camacho2007}. Our motivation to use MPC is that it solves an optimal control problem over a receding time window, which provides the capability of predicting future events and taking actions accordingly. Note that even though this approach is only sub-optimal, in general \cite{Mattingley2010}, it works very well in many applications; our numerical experiments suggest that MPC yields very good performance in traffic control applications. Two schemes are developed in this paper: i) centralized; and ii) decentralized. In the centralized scheme, assuming that the control inputs are inlet traffic flows and configuration of the traffic lights across the network, a two-step control scheme is proposed. In a normal traffic mode, the proposed centralized scheme alleviates traffic density in all lanes, ensuring that traffic density in the entire network is less than a certain value. When an emergency vehicle approaches the network\textemdash this condition is referred as an emergency traffic mode\textemdash the control objective is to clear the path for the emergency vehicle, without causing congestion in other lanes. It is shown that our proposed centralized scheme provides the optimal solution, even though its computation time may be large for large networks. In the decentralized scheme, inlet traffic flows and configuration of the traffic lights at each intersection are controlled by a local control unit, while the control units share data with each other through an aggregator. In the decentralized scheme, the aggregator should receive and send the data during every traffic light state (i.e., `red' or `green'). Since the traffic cycle ranges from one minute to three minutes in real-world traffic systems \cite{NACTO}, the smallest duration of traffic light states is 30 seconds; thus, the maximum allowable communication delay is around 30 seconds, which is achievable even with cheap communication technologies. Thus, the decentralized scheme is more suitable for large networks, even though it yields a sub-optimal solution. Note that the robustness and tolerance of the decentralized scheme to uncertainty in communication delay and communication failures are out of the scope of this paper, and will be considered as future work.

The key contributions of this paper are: i) we develop a traffic control framework which provides an efficient and timely emergency vehicle passage through multiple intersections, without causing congestion in other lanes; ii) we propose a centralized scheme for small networks and a decentralized scheme for large networks that addresses scalability issues in integrating CTM and MPC; and iii) we validate our schemes via extensive simulation studies, and assess their performance in different scenarios. The main features of the proposed framework are: i) it is general and can be applied to any network of interconnected signaled lanes; and ii) it does not require vehicle to everything (V2X) connectivity, and hence it can be utilized in the currently existing traffic systems; the only communication requirement is between the emergency vehicle and the central control unit in the centralized scheme, and with the aggregator in the decentralized scheme. Note that this paper considers only macroscopic characteristics of traffic flow; it is evident that the existence of V2X connectivity can not only be exploited to further improve efficiency at the macro-level, but it can also be leveraged to ensure safety and avoid collisions.

The key innovations of this paper with respect to prior work are: i) formulating the traffic density control problem in both normal and emergency modes as MPC problems; ii) developing a two-step optimization procedure implementable in the current traffic systems; and iii) deriving centralized and decentralized schemes for traffic networks with different size and communication capacity.

The rest of the paper is organized as follows. Section \ref{sec:PS} describes macroscopic discrete-time model of the traffic flow in the network. Section \ref{sec:centralized} discusses the design procedure of the centralized traffic control scheme. The decentralized scheme is discussed in Section \ref{sec:distributed}. Section \ref{sec:simulation} reports simulations results and compares the centralized and decentralized schemes. Finally, Section \ref{sec:conclusion} concludes the paper.

\paragraph*{Notation} $\mathbb{R}$ denotes the set of real numbers, $\mathbb{R}_{\geq0}$ denotes the set of non-negative real numbers, $\mathbb{Z}$ denotes the set of integer numbers, and $\mathbb{Z}_{\geq0}$ denotes the set of non-negative integer numbers. For the matrix $X$, $X^\top$ denotes its transpose, $\rho(X)$ denotes its spectral radius, and $\left\Vert X\right\Vert_1=\sup_{y\neq0}\frac{\left\Vert Xy\right\Vert_1}{\left\Vert y\right\Vert_1}$ with $\left\Vert\cdot\right\Vert_1$ as the $\ell_1$-norm. For the vector $y$, $[y]_+$ is the element-wise rounding to the closest non-negative integer function. For given sets $X,Y$, $X\oplus Y:\{x+y:x\in X,y\in Y\}$ is the Minkowski set sum. TABLE \ref{table:notation} lists the essential notation of this paper.

\begin{table}[!t]
\centering
\caption{List of calligraphic, Greek, Latin, and subscript and superscript symbols.}
\begin{tabular}{l|c|l}
Type & Symbol & Description \\
\hline
Calligraphic & $\mathcal{N}$ & Set of lanes \\
 & $\mathcal{M}$ & Set of intersections \\
 & $\mathcal{G}$ & Traffic network\\
 & $\mathcal{E}$ & Edge of traffic graph\\
 & $\mathcal{X}$ & Constraint set\\
 & $\mathcal{D}$ & Disturbance set\\
\hline
Greek & $\Lambda$ & Set of possible commands by traffic lights\\
& $\lambda$ & Configuration of traffic lights\\
& $\gamma$ & Prioritizing parameter \\
\hline
Latin & $x$ & Traffic density\\
& $y$ & Traffic inflow\\
& $z$ & Traffic outflow\\
& $u$ & Inlet flow \\
& $d$ & Disturbance input\\
& $t$ & Discrete time instant \\
& $k$ & Prediction time instant \\
\hline
Subscript and & $in$ & inlet \\
Superscript & $df$ & Disturbance-free \\
 & $nom$ & Nominal \\
 & $n$ & normal condition \\
 & $e$ & Emergency condition \\
 & $t$ & Computed at time $t$ \\
\hline
\end{tabular}
\label{table:notation}
\end{table}

\section{System Modelling}\label{sec:PS}
In this section, we formulate the traffic control problem for a general traffic network.

\subsection{Traffic Network}\label{sec:graph}
Consider a traffic network with $N$ lanes and $M$ intersections. There are $N_{in}<N$ inlets through which vehicles enter the network. We denote the set of lanes by $\mathcal{N}=\{1,\cdots,N\}$, the set of intersections by $\mathcal{M}=\{1,\cdots,M\}$, and the set of inlets by $\mathcal{N}_{in}\subset\mathcal{N}$.

The considered traffic network can be represented by a graph $\mathcal{G}(\mathcal{N},\mathcal{E})$, where $\mathcal{E}\subset\mathcal{N}\times\mathcal{N}$ defines the edge of graph. The edge $(i,j)\in\mathcal{E}$ represents a directed connection from lane $i$ to lane $j$. Since all lanes are assumed to be unidirectional (note that two-way roads are modeled as two opposite-directional lanes), if $(i,j)\in\mathcal{E}$, we have $(j,i)\not\in\mathcal{E}$. Also, we assume that U-turns are not allowed, i.e., $(i,j),(j,i)\not\in\mathcal{E}$, if lanes $i$ and $j$ are opposite-directional lanes on a single road.

Note that we assume that the traffic graph $\mathcal{G}(\mathcal{N},\mathcal{E})$ remains unchanged; that is we do not consider graph changes due to unexpected events (e.g., changes in the edge $\mathcal{E}$ as a result of lane blockages due to accidents). We leave the developments of strategies for rerouting in the case of a change in the traffic graph to future work.

\subsection{Action Space By Traffic Lights}\label{sec:ActionSpace}
Suppose that all lanes, except outlets, are controlled by traffic lights which have three states: `red', `yellow', and  `green'. The vehicles are allowed to move when the light is `yellow' or `green', while they have to stop when the light is `red'. This means that there are practically two states for each traffic light.

Let $\lambda_j(t)$ be the configuration of traffic lights at intersection $j\in\mathcal{M}$ at time $t$. We denote the set of all possible configurations at intersection $j$ by $\Lambda_j=\{\lambda_{j,1},\cdots,\lambda_{j,\mu_j}\}$, where $\mu_j\in\mathbb{Z}_{\geq0}$. Indeed, the set $\Lambda_j$ represents the set of all possible actions that can be commanded by the traffic lights at intersection $j$. Therefore, the set of all possible actions by traffic lights across the network is $\Lambda=\Lambda_1\times\cdots\times\Lambda_M$, and the $M$-tuple  $\lambda(t)=\big(\lambda_{1}(t),\cdots,\lambda_{M}(t)\big)\in\Lambda$ indicates the action across the network at time $t$.

\subsection{Macroscopic Traffic Flow Model}
The traffic density in each lane is a macroscopic characteristic of traffic flow \cite{Chanut2003,Khan2018}, which can be described by the CTM that transforms the partial differential equations of the macroscopic Lighthill-Whitham-Richards (LWR) model \cite{Yu2021} into simpler difference equations at the cell level. The CTM formulates the relationship between the key traffic flow parameters, and can be cast in a discrete-time state-space form.

Let traffic density be defined as the total number of vehicles in a lane at any time instant, then the traffic inflow is defined as the total number of vehicles entering a lane during a given time period, and traffic outflow is defined as the total number of vehicles leaving a lane during a given time period. We use $x_i(t)\in\mathbb{Z}_{\geq0}$, $y_i(t)\in\mathbb{R}_{\geq0}$, and $z_i(t)\in\mathbb{R}_{\geq0}$ to denote the traffic density, traffic inflow, and traffic outflow in lane $i$ at time $t$, respectively. The traffic dynamics \cite{Adacher2018,Vishnoi2020} in lane $i$ can be expressed as
\begin{align}\label{eq:LWRlanei}
x_i(t+1)=\left[x_i(t)+y_i(t)-z_i(t)\right]_+,
\end{align}
where the time interval $[t,t+1)$ is equivalent to $\Delta T$ seconds. Since $x_i(t)$ is defined as the number of existing vehicles in each lane, we use the rounding function in \eqref{eq:LWRlanei} to ensure that $x_i(t)$ remains a non-negative integer at all times. Given $\Delta T$,  $y_i(t)$ and $z_i(t)$ are equal to the number of vehicles entering and leaving the lane $i$ in $\Delta T$ seconds, respectively. 

The traffic outflow $z_i(t)$ can be computed as \cite{Rastgoftar2021}
\begin{align}\label{eq:zi}
z_i(t)=p_i\big(\lambda(t)\big)x_i(t),
\end{align}
where $p_i\big(\lambda(t)\big)$ is the fraction of outflow vehicles in lane $i$ during the time interval $[t,t+1)$, satisfying 
\begin{align}\label{eq:pi}
p_i\big(\lambda(t)\big)\left\{
\begin{array}{ll}
   =0,  & \text{if traffic light of lane }i\text{ is `red'} \\
  \in[0,1],  & \text{if traffic light of lane }i\text{ is `green'}
\end{array}
\right.;
\end{align}
in other words, $p_i\big(\lambda(t)\big)$ is the ratio of vehicles leaving lane $i$ during the time interval $[t,t+1)$ to the total number of vehicles in lane $i$ at time instant $t$. It is noteworthy that even though the impact of lane blockage or an accident in lane $i$ can be modeled by adjusting $p_i\big(\lambda(t)\big)$, this paper does not aim to deal with such unexpected events.

\begin{remark}
We assume that outlet traffic flows are uncontrolled, i.e., there is no traffic light or gate at the end of outlets. This assumption is plausible, as any road connecting the considered traffic network to the rest of the grid can be divided at a macro-level into an uncontrollable outlet inside the considered network and a lane outside the considered network (possibly controlled with a traffic light or a network gate). The extension of the proposed methods to deal with controlled outlet flows is straightforward by modifying \eqref{eq:zi} and all presented optimization problems to account for outlet flow (similar to what we do for inlet flow $u_i(t)$); thus, to simplify the exposition and subsequent developments, we will not discuss controlled outlets.
\end{remark}

The traffic inflow $y_i(t)$ can be computed as
\begin{align}\label{eq:yi}
y_i(t)=\left\{
\begin{array}{ll}
   u_i(t),  & \text{if $i\in\mathcal{N}_{in}$} \\
   \sum\limits_{j=1}^Nq_{j,i}\big(\lambda(t)\big)z_j(t),  & \text{otherwise}
\end{array}
\right.,
\end{align}
where $u_i(t)\in\mathbb{Z}_{\geq0}$ is the inlet flow which is defined as the number of vehicles entering the traffic network through inlet $i$ during the time interval $[t,t+1)$. The computed optimal inflows can be implemented by means of network gates, i.e., ramp meters \cite{Gomez2006,Gomez2008} for highways and metering gates \cite{Mohebifard2018} for urban streets). In \eqref{eq:yi}, $q_{j,i}\big(\lambda(t)\big)$ is the fraction of outflow of lane $j$ directed toward lane $i$ during the time interval $[t,t+1)$, which is
\begin{align}\label{eq:qji}
q_{j,i}\big(\lambda(t)\big)\left\{
\begin{array}{ll}
   =0,  & \text{(if traffic light of lane $i$ is `red')} \\
   & \text{OR (if }(j,i)\not\in\mathcal{E}\text{)}\\
   \in[0,1],  & \text{(if traffic light of lane $i$ is `green')}\\
   & \text{AND (if }(j,i)\in\mathcal{E}\text{)}
\end{array}
\right.,
\end{align}
and satisfies $\sum_{i=1}^Nq_{j,i}\big(\lambda(t)\big)=1$ for all $j\in\mathcal{N}$. More precisely, $q_{j,i}\big(\lambda(t)\big)$ is the ratio of vehicles leaving lane $j$ and entering lane $i$ during the time interval $[t,t+1)$ to the total number of vehicles leaving lane $j$ during the time interval $[t,t+1)$.

From \eqref{eq:LWRlanei}-\eqref{eq:qji}, traffic dynamics of the entire network can be expressed as
\begin{align}\label{eq:system1}
x(t+1)=\left[\bar{A}\big(\lambda(t)\big)x(t)+B\bar{U}(t)\right]_+,
\end{align}
where $x(t)=[x_1(t)~\cdots~x_N(t)]^\top\in\mathbb{Z}_{\geq0}^N$, $\bar{A}:\Lambda\rightarrow\mathbb{R}^{N\times N}$ is the so-called \textit{traffic tendency matrix} \cite{Rastgoftar2019}, $B\in\mathbb{R}^{N\times N_{in}}$, and $\bar{U}(t)\in\mathbb{Z}_{\geq0}^{N_{in}}$ is the boundary inflow vector. It should be noted that the $(i,j)$ element of $B$ is 1 if lane $i$ is the $j$-th inlet, and 0 otherwise.

\begin{remark}
At any $t$, the $(i,i)$ element of the traffic tendency matrix $\bar{A}\big(\lambda(t)\big)$ is $1-p_i\big(\lambda(t)\big)$. Also, its $(i,j)$ element ($i\neq j$) is $q_{j,i}\big(\lambda(t)\big)p_j\big(\lambda(t)\big)$. As a result, since $\sum_{i=1}^Nq_{j,i}\big(\lambda(t)\big)=1,~\forall j\in\mathcal{N}$, the maximum absolute column sum of the traffic tendency matrix is less than or equal to 1. This means that at any $t$, we have $\left\Vert\bar{A}\big(\lambda(t)\big)\right\Vert_1\leq1$, which implies that $\rho\Big(\bar{A}\big(\lambda(t)\big)\Big)\leq1$. Therefore, $\rho\Big(\bar{A}\big(\lambda(t)\big)\bar{A}\big(\lambda(t+1)\big)\bar{A}\big(\lambda(t+2)\big)\cdots\Big)\leq1$, which means that the unforced system (i.e., when $\bar{U}(t)=0$) is stable, although trajectories may not asymptotically converge to the origin. This conclusion is consistent with the observation that in the absence of new vehicles entering to lane $i$, the traffic density in lane $i$ remains unchanged if the corresponding traffic light remains `red'. 
\end{remark}

\begin{remark}
In general, system \eqref{eq:system1} is not bounded-input-bounded-output stable. For instance, the traffic density in lane $i$ constantly increases if $y_i(t)>0$ at all times and the corresponding traffic light remains `red'. 
\end{remark}

Given the action $\lambda(t)$, the traffic dynamics given in \eqref{eq:system1} depend on the parameters $p_i\big(\lambda(t)\big)$ and $q_{j,i}\big(\lambda(t)\big),~\forall i,j$, as well as the boundary inflow vector $\bar{U}(t)$. These parameters are, in general, \textit{a priori} unknown. We assume that these parameters belong to some bounded intervals, and we can estimate these intervals from prior traffic data. Thus, traffic dynamics given in \eqref{eq:system1} can be rewritten as
\begin{align}\label{eq:system2}
x(t+1)=\Big[\Big(A\big(\lambda(t)\big)&+\Delta A(t)\Big)x(t)\nonumber\\
&+B\big(U(t)+\Delta U(t)\big)\Big]_+,
\end{align}
where $A\big(\lambda(t)\big)\in\mathbb{R}^{N\times N}$ is the traffic tendency matrix computed by nominal values of $p_i$ and $q_{j,i},~\forall i,j$ associated with the action $\lambda(t)$, $\Delta A(t)\in\mathbb{R}^{N\times N}$ covers possible uncertainties, $U(t)\in\mathbb{Z}_{\geq0}^{N_{in}}$ is the boundary inflow vector at time $t$, and $\Delta U(t)\in\mathbb{Z}_{\geq0}^{N_{in}}$ models possible inflow uncertainties.

\begin{remark}
The boundary inflow $U(t)$ is either uncontrolled or controlled. In the case of uncontrolled inlets, $U(t)$ represents the nominal inflow learnt from prior data, which, in general, is time-dependent, as it can be learnt for different time intervals in a day (e.g., in the morning, in the evening, etc). In this case, $\Delta U(t)$ models possible imperfections. In the case of a controlled inlet traffic flows, $U(t)$ is the control input at time $t$. Note that $U(t)$ determines the available throughput in inlets, i.e., an upper-bound on vehicles entering the network through each inlet. However, traffic demand might be less than the computed upper-bounds, meaning that utilized throughput is less than the available throughput. In this case, $\Delta U(t)$ models differences between the available and utilized throughput. 
\end{remark}

Finally, due to the rounding function in \eqref{eq:system2}, the impact of the uncertainty terms $\Delta A(t)$ and $\Delta U(t)$ can be expressed as an additive integer. More precisely, traffic dynamics given in \eqref{eq:system2} can be rewritten as
\begin{align}\label{eq:system3}
x(t+1)=\max\Big\{&\Big[A\big(\lambda(t)\big)x(t)+BU(t)\Big]_++d(t),0\Big\},
\end{align}
where $d(t)=[d_1(t)~\cdots~d_N(t)]^\top\in\mathcal{D},~\forall t$ is the disturbance that is unknown but bounded, with $\mathcal{D}\subset\mathbb{Z}^{N}$ as a polyhedron containing the origin. Note that $d_i(t)$ also models vehicles parking/unparking in lane $i$.

\section{Emergency Vehicle-Centered Traffic Control\textemdash Centralized Scheme}\label{sec:centralized}
\begin{figure}[!t]
\centering
\includegraphics[width=8cm]{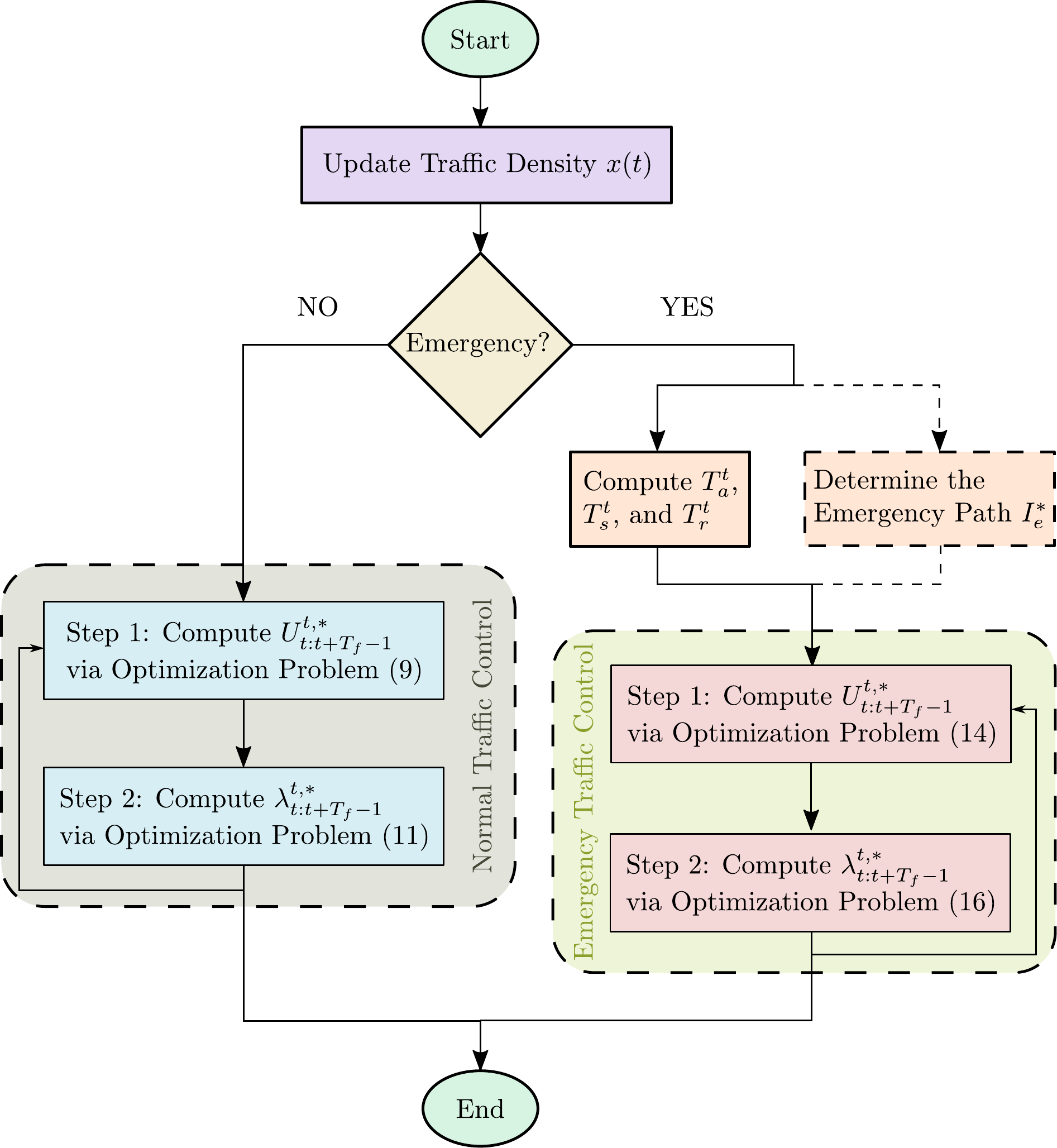}
\caption{Algorithmic flowchart of the proposed centralized traffic control scheme. This algorithm should be run at every $t$ in the central control unit.}
\label{fig:ControlSchemeCentralized}
\end{figure}

In this section, we will propose a centralized scheme whose algorithmic flowchart given in Fig. \ref{fig:ControlSchemeCentralized}. As seen in this figure, a central control unit determines the optimal inlet flows and configuration of all traffic lights. This implies that the data from all over the network should be available to the central unit at any $t$.

In this section, we will use the following notations. Given the prediction horizon $[t,t+T_f]$ for some $T_f\in\mathbb{Z}_{\geq0}$, we define $U_{t:t+T_f-1}^{t}=[U^{t}(t)^\top~\cdots~U^{t}(t+T_f-1)^\top]^\top\in\mathbb{Z}_{\geq0}^{T_fN_{in}}$, where $U^{t}(t+k)\in\mathbb{Z}_{\geq0}^{N_{in}}$ is the boundary inflow vector for time $t+k$ (with $k\leq T_f-1$) computed at time $t$. Also, $\lambda_{t:t+T_f-1}^{t}=\{\lambda^{t}(t),\cdots,\lambda^{t}(t+T_f-1)\}\in\Lambda^{T_f}$, where $\lambda^{t}(t+k)$ is the configuration of all traffic lights for time $t+k$ (with $k\leq T_f-1$) computed at time $t$. Note that $\ast$ is added to the above-mentioned notations to indicate optimal decisions.

\subsection{Normal Traffic Mode}\label{sec:NormalCondition}
The normal traffic mode corresponds to traffic scenarios in which there is no emergency vehicle. Given the prediction horizon $[t,t+T_f]$, the control objective in a normal traffic mode is to determine boundary inflows and configurations of traffic lights over the prediction horizon such that traffic congestion is alleviated in all lanes. This objective can be achieved through the following two-step receding horizon control; that is, the central unit computes the optimal boundary inflows and configuration of traffic lights over the prediction horizon by solving the associated optimization problems at every time instant $t$, but only implements the next boundary inflows and configuration of traffic lights, and then solves the associated optimization problems again at the next time instant, repeatedly.

\subsubsection{Step 1}
Consider $\{\lambda_{t:t+T_f-2}^{t-1,\ast},\lambda(t+T_f-1)\}$, where $\lambda_{t:t+T_f-2}^{t-1,\ast}$ is the optimal solution\footnote{$\lambda_{0:T_f-2}^{-1,\ast}$ should be selected randomly from the action space $\Lambda^{T_f-1}$.} of \eqref{eq:OptStep2} obtained at time $t-1$ and $\lambda(t+T_f-1)$ is selected randomly from the action space $\Lambda$. Then, the optimal boundary inflows over the prediction horizon $[t,t+T_f]$ (i.e., $U_{t:t+T_f-1}^{t,\ast}$) can be obtained by solving the following optimization problem:
\begin{subequations}\label{eq:OptStep1}
\begin{align}\label{eq:OptStep1_a}
\min\limits_{U}\; \sum\limits_{k=0}^{T_f-1}\Big(\left\Vert \hat{x}_{df}(k|t)\right\Vert_{\Gamma_n}^2+\left\Vert U(t+k)-U_{nom}(t+k)\right\Vert_{\Theta}^2\Big),
\end{align}
subject to
\begin{align}
\hat{x}(k|t)\subseteq\hat{\mathcal{X}},&~k=1,\cdots,T_f,\label{eq:OptStep1_b}\\
U(t+k)\in\mathbb{Z}_{\geq0}^{N_{in}},&~k=0,\cdots,T_f-1,\label{eq:OptStep1_c}
\end{align}
\end{subequations}
where $\Theta=\Theta^\top\geq0$ ($\in\mathbb{R}^{N_{in}\times N_{in}}$) is a weighting matrix, $\hat{\mathcal{X}}\subset\mathbb{Z}_{\geq0}^N$ is a polyhedron containing the origin\footnote{The upper-bound on the traffic density of each lane can be specified according to the capacity of the lane. See \cite{Makki2020} for
a comprehensive survey.}, and
\begin{align}\label{eq:SystemApproximation}
\hat{x}(k+1|t)\in\Big(A\big(\lambda^{t-1,\ast}(t+k)\big)\hat{x}(k|t)+BU(t+k)\Big)\oplus\mathcal{D},
\end{align}
with initial condition $\hat{x}(0|t)=x(t)$, and $\lambda^{t-1,\ast}(t+T_f-1)=\lambda(t+T_f-1)$ which is selected randomly from the action space $\Lambda$. Note that to account for the disturbance $d(t)$, \eqref{eq:SystemApproximation} uses the Minkowski set-sum of nominal predictions plus the set of all possible effects of the disturbance $d(t)$ on the traffic density. The subscript ``df" in \eqref{eq:OptStep1_a} stands for disturbance-free, and $\hat{x}_{df}(k|t)$ can be computed via \eqref{eq:SystemApproximation} by setting $\mathcal{D}=\{\textbf{0}\}$. The $U_{nom}(t)$ is the nominal boundary inflow at time $t$, which can be estimated based on prior traffic data. In \eqref{eq:OptStep1}, $\Gamma_n=\text{diag}\{\gamma_{1}^n,\cdots,\gamma_{N}^n\}$, where $\gamma_i^n\geq0,~\forall i\in\{1,\cdots,N\}$ is a design parameter that can be used to prioritize lanes. As suggested by the U.S. Department of Transportation \cite{DoT_priority}, the prioritizing parameters can be determined according to total crashes and congestion over a specified period of time (e.g., over a 5-year period); the higher the prioritizing parameter is, the higher priority is given to the density alleviation.

In summary, Step 1 computes the optimal boundary inflows by solving the optimization problem \eqref{eq:OptStep1}, which has $T_f\times N_{in}$ integer decision variables constrained to be non-negative, and has $T_f\times N$ inequality constraints on traffic density.

\subsubsection{Step 2}
Given $U_{t:t+T_f-1}^{t,\ast}$ as the optimal solution of \eqref{eq:OptStep1} obtained at time $t$, the optimal configuration of all traffic lights over the prediction horizon $[t,t+T_f]$ (i.e., $\lambda_{t:t+T_f-1}^{t,\ast}$) can be determined by solving the following optimization problem:
\begin{subequations}\label{eq:OptStep2}
\begin{align}
\min\limits_{\lambda}\; \sum\limits_{k=0}^{T_f-1}\left\Vert \tilde{x}_{df}(k|t)\right\Vert_{\Gamma_n}^2,
\end{align}
subject to
\begin{align}
\tilde{x}(k|t)\subseteq\tilde{\mathcal{X}},&~k=1,\cdots,T_f,\\
\lambda(t+k)\in\Lambda,&~k=0,\cdots,T_f-1,
\end{align}
\end{subequations}
where $\tilde{\mathcal{X}}\subset\mathbb{Z}_{\geq0}^N$ is a polyhedron containing the origin, and
\begin{align}\label{eq:xtilde}
\tilde{x}(k+1|t)\in\max\Big\{\Big[A\big(\lambda&(t+k)\big)\tilde{x}(k|t)\nonumber\\
&+BU^{t,\ast}(t+k)\Big]_+\oplus\mathcal{D},0\Big\},
\end{align}
with the initial condition, $\tilde{x}(0|t)=x(t)$. Note that $\tilde{x}_{df}(k|t)$ can be computed via \eqref{eq:xtilde} by setting $\mathcal{D}=\{\textbf{0}\}$. Note that similar to \eqref{eq:SystemApproximation}, a set-valued prediction of traffic density by taking into account all possible realizations of the disturbance $d(t)$ is considered in \eqref{eq:xtilde} to account for the disturbance $d(t)$.

In summary, Step 2 determines the optimal configuration of traffic lights across the network by solving the optimization problem \eqref{eq:OptStep2} which has $T_f$ decision variables (each one is an $M$-tuple representing the configuration of traffic lights) constrained to belong to the set $\Lambda$ (see Subsection \ref{sec:ActionSpace}), and has $T_f\times N$ inequality constraints on traffic density.

\begin{remark}
The cost function in \eqref{eq:OptStep1} has two terms. The first term penalizes traffic density in all lanes of the network, and the second term penalizes the difference between the inlet traffic flows and their nominal values. It should be noted that a sufficiently large matrix $\Theta$ guarantees that vehicles will never be blocked behind the network gates. A different method \cite{Rastgoftar2021} to ensure that vehicles will not be blocked is to constrain the total boundary inflow to be equal to a certain amount, i.e., $\sum_{i\in\mathcal{N}_{in}} u_i(t)=\bar{u},~\forall t$, where $\bar{u}$ can be determined based upon prior traffic data. It is noteworthy that the computed optimal inflows can be implemented by means of network gates, i.e., ramp meters \cite{Gomez2006,Gomez2008} for highways and metering gates \cite{Mohebifard2018} for urban streets).
\end{remark}

\begin{remark}
The prediction given in \eqref{eq:SystemApproximation} provides an approximation to system \eqref{eq:system3}, and the traffic density may take non-integer and/or negative values. However, as will be shown later, this approximation is efficient in ensuring optimality. The main advantage of using such an approximation is that the integer  programming as in \eqref{eq:OptStep1} can be easily solved by available tools.
\end{remark}

\begin{remark}
The optimization problem \eqref{eq:OptStep2} can be solved by using the brute-force search \cite{Mahoor2017} (a.k.a. exhaustive search or generate\&test) algorithm. Note that the size of the problem \eqref{eq:OptStep2} is limited, since $\Lambda^{T_f}$ and $\mathcal{D}$ are finite. However, there are some techniques to reduce the search space, and consequently speed up the algorithm. For instance, if the configuration $\lambda_{t:t+T_f-1}^{t}$ is infeasible and causes congestion at time $t+k$ (with $0\leq k\leq T_f-1$), all configurations with the same first $k-1$ actions will be excluded from the search space. Our simulation studies show that this simple step can largely reduce the computation time of the optimization problem \eqref{eq:OptStep2} (in our case, from 10 seconds to 8 milliseconds). 
\end{remark}

\begin{remark}
In the case of uncontrolled boundary inflow, the proposed scheme for normal traffic mode reduces to solving only the optimization problem \eqref{eq:OptStep2} based upon learnt nominal boundary inflows. 
\end{remark}

\begin{remark}
We assume that constraints on the traffic density are defined such that the resulting optimization problems are feasible. However, in the case of infeasibility, we can use standard methods (e.g., introducing slack variables) to relax constraints.
\end{remark}

\subsection{Emergency Traffic Mode}\label{sec:EmergencyCondition}
Suppose that:
\begin{itemize}
\item At time $t=t_e$, a notification is received by the central control unit indicating that an emergency vehicle will enter the network in $T_a^t$ time steps. Note that for $t<t_e$ the condition of the network was normal.
\item Given the entering and leaving lanes, let $\mathcal{P}$ represents the set of all possible paths for the emergency vehicle. Once the notification is received, i.e., at time $t=t_e$, based on the current and predicted traffic conditions, the optimal emergency path $I_e^\ast$ should be selected by the central control unit (see Remark \ref{remkr:EmergencyPath}) and be given to the emergency vehicle. We assume that the emergency vehicle will follow the provided path. 
\item The emergency vehicle should leave the network in maximum $T_s^t$ time steps.
\item Once the emergency vehicle leaves the network, the traffic density in all lanes should be recovered to the normal traffic mode in $T_r^t$ time steps. This phase will be referred as the recovery phase in the rest of the paper. 
\end{itemize}

\begin{remark}
$T_a^t$, $T_s^t$, and $T_r^t$ are specified at time $t$. These values can be computed by leveraging connectivity between the emergency vehicle and the roadside infrastructure. Note that these variables are time-variant, as they should be recomputed based on the traffic condition and position of the emergency vehicle at any $t$. For instance, once the emergency vehicle enters the network, $T_a^t$ should be set to zero, and once the emergency vehicle leaves the network $T_s^t$ should be set to zero. Also, when the recovery phase ends, $T_r^t$ will be zero. 
\end{remark}


The control objective in an emergency traffic mode is to shorten the traveling time of the emergency vehicle, i.e., to help the emergency vehicle traverse the network as quickly and efficiently as possible. Given the emergency path with length $L$, the traveling time of the emergency vehicle can be estimated \cite{Zhao2008,Zhang2013} as  
\begin{align}\label{eq:TravelingTime}
\text{Traveling Time}=\frac{L}{V_d}+\beta\times\text{Traffic Density on the Path},
\end{align}
for some constant $\beta>0$, where $V_d$ is the desired traverse velocity. This relationship indicates that for fixed $L$ and $V_d$, to shorten the traveling time of the emergency vehicle one would need to reduce the traffic density on the emergency path.

Therefore, in an emergency traffic mode, given the prediction horizon $[t,t+T_f]$ with $T_f\geq T_a^t+T_s^t+T_r^t$, the control objective can be achieved by determining boundary inflows and configuration of all traffic lights such that: i) during the time interval $[t,t+T_a^t+T_s^t]$ traffic density in emergency path should be reduced as much as possible, while traffic density in other lanes is less than a certain amount; ii) during the time interval $[t+T_a^t+T_s^t,t+T_a^t+T_s^t+T_r^t]$ the traffic density in all lanes should be recovered to the normal traffic mode; and iii) during the time interval $[t+T_a^t+T_s^t+T_r^t,T_f]$ the traffic density in all lanes should satisfy constraints of normal mode.

We propose the following two-step receding horizon control approach to satisfy the above-mentioned objectives. In this approach, the central unit computes the optimal boundary inflows and configuration of traffic lights over the prediction horizon by solving the associated optimization problems at every time instant $t$, but only implements the next boundary inflows and configuration of traffic lights, and then solves the associated optimization problems again at the next time instant, repeatedly.

\subsubsection{Step 1} Consider $\{\lambda_{t:t+T_f-2}^{t-1,\ast},\lambda(t+T_f-1)\}$, where $\lambda_{t:t+T_f-2}^{t-1,\ast}$ is the optimal solution\footnote{Since the traffic condition was normal for $t<t_e$, $\lambda_{0:T_f-2}^{t_e-1,\ast}$ is the optimal solution of \eqref{eq:OptStep2} at time $t_e-1$.} of \eqref{eq:OptStep4} obtained at time $t-1$ and $\lambda(t+T_f-1)$ is selected randomly from the action space $\Lambda$. Then, the optimal boundary inflows over the prediction horizon $[t,t+T_f]$ (i.e., $U_{t:t+T_f-1}^{t,\ast}$) can be computed by solving the following optimization problem:
\begin{subequations}\label{eq:OptStep3}
\begin{align}
\min\limits_{U}\; \sum\limits_{k=0}^{T_f-1}\Big(\left\Vert \hat{x}_{df}(k|t)\right\Vert_{\Gamma_e}^2+\left\Vert U(t+k)-U_{nom}(t+k)\right\Vert_{\Theta}^2\Big),
\end{align}
subject to
\begin{align}
\hat{x}(k|t)\subseteq\hat{\mathcal{X}}^+,&~k=1,\cdots,T_a^t+T_s^t+T_r^t\\
\hat{x}(k|t)\subseteq\hat{\mathcal{X}},&~k=T_a^t+T_s^t+T_r^t+1,\cdots,T_f\\
U(t+k)\in\mathbb{Z}_{\geq0}^{N_{in}},&~k=0,\cdots,T_f-1,
\end{align}
\end{subequations}
where $\hat{x}(k|t)$ is as in \eqref{eq:SystemApproximation},  $\hat{\mathcal{X}}^+\supset\hat{\mathcal{X}}$ is the extended constraint set (see Remark \ref{remark:Extension}), and $\Gamma_e=\text{diag}\{\gamma_1^e(k),\cdots,\gamma_N^e(k)\}$ (see Remark \ref{remark:Gammae}) with
\begin{align}\label{eq:Gammae}
\gamma_i^e(k)=\left\{
\begin{array}{ll}
    \bar{\gamma}_e, & \text{if $i\in I_e^\ast$ and $k\leq T_a^t+T_s^t$} \\
   \gamma_i^n,  & \text{otherwise}
\end{array}
\right.,
\end{align}
with $\bar{\gamma}_e\gg\max_{i}\{\gamma_i^n\}$, and $I_e^\ast$ is the selected emergency path (see Remark \ref{remkr:EmergencyPath}). The prioritizing parameters as in \eqref{eq:Gammae} ensure that the traffic density in the lanes included in the emergency path will be alleviated with a higher priority in the emergency traffic mode.

Similar to \eqref{eq:OptStep1}, the optimization problem \eqref{eq:OptStep3} has $T_f\times N_{in}$ integer decision variables constrained to be non-negative, and has $T_f\times N$ inequality constraints on traffic density.

\subsubsection{Step 2}
Given $U_{t:t+T_f-1}^{t,\ast}$ as the optimal solution of \eqref{eq:OptStep3} obtained at time $t$, the optimal configurations of the traffic lights over the prediction horizon $[t,t+T_f]$ (i.e., $\lambda_{t:t+T_f-1}^{t,\ast}$) can be determined by solving the following optimization problem:
\begin{subequations}\label{eq:OptStep4}
\begin{align}
\min\limits_{\lambda}\; \sum\limits_{k=0}^{T_f-1}\left\Vert \tilde{x}_{df}(k|t)\right\Vert_{\Gamma_e}^2,
\end{align}
subject to
\begin{align}
\tilde{x}(k|t)\subseteq\tilde{\mathcal{X}}^+,&~k=1,\cdots,T_a^t+T_s^t+T_r^t,\\
\tilde{x}(k|t)\subseteq\tilde{\mathcal{X}},&~k=T_a^t+T_s^t+T_r^t+1,\cdots,T_f\\
\lambda(t+k)\in\Lambda,&~k=0,\cdots,T_f-1
\end{align}
\end{subequations}
where $\tilde{x}(k|t)$ is as in \eqref{eq:xtilde}, and $\tilde{\mathcal{X}}^+\supset\tilde{\mathcal{X}}$ is the extended set (see Remark \ref{remark:Extension}). Similar to \eqref{eq:OptStep2}, the optimization problem \eqref{eq:OptStep4} has $T_f$ decision variables (each one is an $M$-tuple representing the configuration of traffic lights) constrained to belong to the set $\Lambda$ (see Subsection \ref{sec:ActionSpace}), and has $T_f\times N$ inequality constraints on traffic density.

\begin{remark}
The optimization problem \eqref{eq:OptStep3} can be solved by mixed-integer tools, and the optimization problem \eqref{eq:OptStep4} can be solved by using the brute-force search algorithms.
\end{remark}

\begin{remark}\label{remark:Extension}
We assume that constraints on the traffic density can be temporarily relaxed. This assumption is reasonable \cite{Kolmanovsky2011,Li2021}, as in practice, constraints are often imposed conservatively to avoid congestion. In mathematical terms, by relaxation we mean that 
traffic density should belong to extended sets $\hat{\mathcal{X}}^+\supset\hat{\mathcal{X}}$ and $\tilde{\mathcal{X}}^+\supset\tilde{\mathcal{X}}$. This relaxation enables the control scheme to put more efforts on alleviation of traffic density in emergency path. This relaxation can last up to maximum $T_a^{t_e}+T_s^{t_e}+T_r^{t_e}$ time steps. 
\end{remark}

\begin{remark}\label{remark:Gammae}
$\Gamma_e$ as in \eqref{eq:Gammae} prioritizes alleviating traffic density in lanes included in the emergency path $I_e^\ast$ during the time interval in which the emergency vehicle is traversing the network, i.e., the time interval $[t,t+T_a^t+T_s^t]$.
\end{remark}

\begin{remark}\label{remkr:EmergencyPath}
Once the emergency notification is received by the central control unit (i.e., at time $t=t_e$), the optimization problems \eqref{eq:OptStep3} and \eqref{eq:OptStep4} should be solved for all possible paths, i.e., for each element of $\mathcal{P}$. Then: i) according to \eqref{eq:TravelingTime}, the optimal emergency path $I_e^\ast$ should be selected as
\begin{align}
I_e^\ast=\arg\;\min\limits_{I_e\in\mathcal{P}}\;\sum\limits_{k=1}^{T_a^{t_e}+T_s^{t_e}}\sum\limits_{i\in I_e}x_{i,df}(k|t_e);
\end{align}
and ii) the boundary inflow and configuration of traffic lights at time $t=t_e$ will be the ones associated with the optimal emergency path $I_e^\ast$. 
\end{remark}

\begin{remark}
Once the recovery phase ends, the traffic condition will be normal, and the boundary inflow vector and configuration of traffic lights should be determined through the two-step control scheme presented in Subsection \ref{sec:NormalCondition}
\end{remark}

\begin{remark}
In the case of uncontrolled boundary inflow, the proposed scheme for emergency traffic mode reduces to solving only the optimization problem \eqref{eq:OptStep4} based upon learnt nominal boundary inflows.
\end{remark}

\section{Emergency Vehicle-Centered Traffic Control\textemdash Decentralized Scheme}\label{sec:distributed}
\begin{figure}[!t]
\centering
\includegraphics[width=8cm]{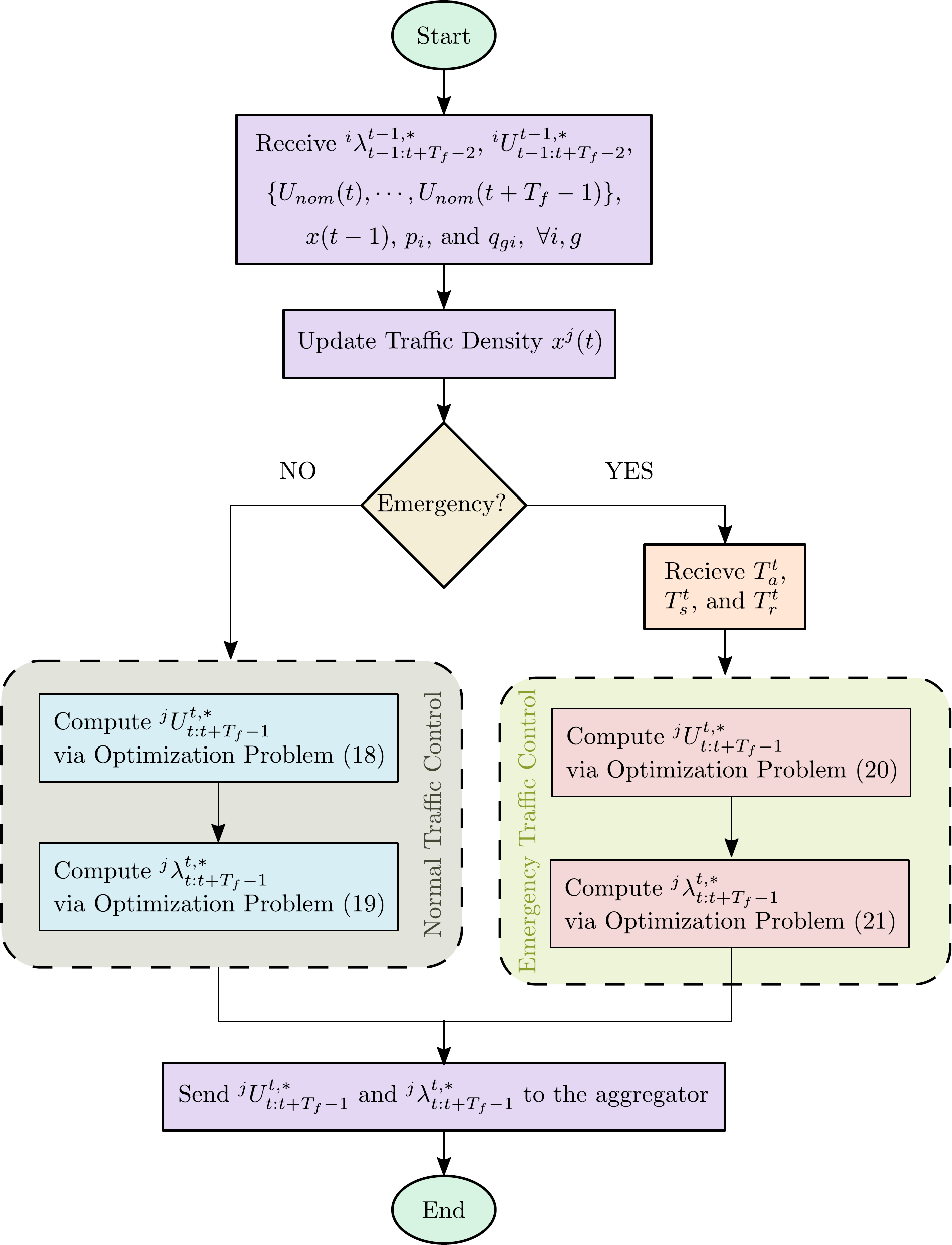}
\caption{Algorithmic flowchart of the proposed decentralized traffic control scheme. This algorithm should be run at every $t$ in the Control Unit\#$j$.}
\label{fig:ControlSchemeDistributed}
\end{figure}

In this section, we will develop a decentralized traffic control scheme whose algorithmic flowchart is depicted in Fig. \ref{fig:ControlSchemeDistributed}. In decentralized scheme, there is a control unit at each intersection, which controls  configuration of the traffic lights at that intersection, as well as the traffic flow in the corresponding inlets. During each sampling period, an aggregator receives data from all control units, augments data, and shares across the network. This is reasonable for real-time applications even with cheap and relatively high-latency communication technologies, as the duration of the traffic light states is large (e.g., 30 seconds). In Section \ref{sec:simulation}, we will characterize the optimality of the developed decentralized scheme in our numerical experiments in different traffic modes in comparison with the centralized scheme.

The main advantage of the decentralized scheme is that the size of resulting optimization problem is very small compared to that of centralized scheme, as it only needs to determine the configuration of traffic lights and inlet traffic flows at one intersection. This greatly reduces the computation time for large networks, even though it may slightly degrade performance. This will be discussed in Section \ref{sec:simulation}.

In this section, we use $^jx(t)\in\mathbb{R}^{N_j},~j\in\mathcal{M}$ (with $N_j\leq N$) to denote traffic density in lanes controlled by Control Unit\#$j$. Also, $^jU_{t:t+T_f-1}^{t}=[^jU^{t}(t)^\top~\cdots~^jU^{t}(t+T_f-1)^\top]^\top\in\mathbb{Z}_{\geq0}^{T_fN^j_{in}},~j\in\mathcal{M}$, where $^jU^{t}(t+k)\in\mathbb{Z}_{\geq0}^{N^j_{in}}$ is the traffic flows in inlets associated with intersection $I_j$ for time $t+k$ (with $k\leq T_f-1$) computed at time $t$, and $^j\lambda_{t:t+T_f-1}^{t}=\{\lambda_j^{t}(t),\cdots,\lambda_j^{t}(t+T_f-1)\}\in\Lambda_j^{T_f},~j\in\mathcal{M}$, where $\lambda_j^{t}(t+k)$ is the configuration of traffic lights at intersection $I_j$ for time $t+k$ (with $k\leq T_f-1$) computed at time $t$. Note that $\sum_jN_{in}^j=N_{in}$, and $\ast$ in the superscript of the above-mentioned notations indicates optimal decisions.

\subsection{Normal Traffic Mode}
As discussed in Subsection \ref{sec:NormalCondition}, the control objective in a normal traffic mode is to alleviate traffic density across the network. During the time interval $[t-1,t)$, all control units receive $^i\lambda_{t-1:t+T_f-2}^{t-1,\ast}$ and $^iU_{t-1:t+T_f-2}^{t-1,\ast}$ for all $i\in\mathcal{M}$, $x(t-1)$, $\{U_{nom}(t),\cdots,U_{nom}(t+T_f-1)\}$, and $p_i$ and $q_{g,i},~i,g\in\mathcal{N}$ from the aggregator. At any $t$, the Control Unit\#$j,~j\in\mathcal{M}$ follows the following steps to determine the inlet traffic flows and the configuration of the traffic lights at intersection $I_j$ in a normal traffic mode:
\begin{enumerate}
\item Compute $x(t|t-1)$ based on the shared information by the aggregator, and according to \eqref{eq:system3} with $d(t-1)=0$.
\item Update traffic density at local lanes (i.e., $^jx(t)$), and replace corresponding elements in $x(t|t-1)$ with updated values.
\item Compute $\{^i\lambda_{t:t+T_f-2}^{t-1,\ast},\lambda_i(t+T_f-1)\}$ for all $i\in\mathcal{M}$, where $^i\lambda_{t:t+T_f-2}^{t-1,\ast}$ is the optimal solution\footnote{$^i\lambda_{0:T_f-2}^{-1,\ast}$ should be selected randomly from the action space $\Lambda_i^{T_f-1}$.} of Control Unit\#$i$ obtained at time $t-1$ and $\lambda_i(t+T_f-1)$ is selected randomly from the action space $\Lambda_i$. 
\item Compute $\{^iU_{t:t+T_f-2}^{t-1,\ast},^iU_{nom}(t+T_f-1)\}$ for all $i\in\mathcal{M}$ and $i\neq j$, where $^iU_{t:t+T_f-2}^{t-1,\ast}$ is the optimal solution\footnote{$^iU_{0:T_f-2}^{-1,\ast}$ is $\{^iU_{nom}(0),\cdots,^iU_{nom}(T_f-2)\}$.} of Control Unit\#$i$ obtained at time $t-1$. 
\item Solve the following optimization problem to determine the inlet traffic flows at intersection $I_j$ over the prediction horizon $[t,t+T_f]$ (i.e., $^jU_{t:t+T_f-1}^{t,\ast}$):
\begin{subequations}\label{eq:OptStep5}
\begin{align}
\min\limits_{^jU}\; \sum\limits_{k=0}^{T_f-1}\Big(&\left\Vert ^j\hat{x}_{df}(k|t)\right\Vert_{^j\Gamma_n}^2\nonumber\\
&+\left\Vert ^jU(t+k)-~^jU_{nom}(t+k)\right\Vert_{^j\Theta}^2\Big),
\end{align}
subject to
\begin{align}
^j\hat{x}(k|t)\subseteq~^j\hat{\mathcal{X}},&~k=1,\cdots,T_f,\\
^jU(t+k)\in\mathbb{Z}_{\geq0}^{N_{in}^j},&~k=0,\cdots,T_f-1,
\end{align}
\end{subequations}
where $^j\Gamma_n=~^j\Gamma_n^\top\geq0$ ($\in\mathbb{R}^{N_j}$) and $^j\Theta=~^j\Theta^\top\geq0$ ($\in\mathbb{R}^{N_{in}^j\times N_{in}^j}$)
are weighting matrices, 
$^j\hat{x}(k|t)$ can be computed via \eqref{eq:SystemApproximation} with initial condition $x(t|t-1)$, and  $^j\hat{\mathcal{X}}\subset\mathbb{R}_{\geq0}^{N_j}$ is a polyhedron containing the origin. The optimization problem \eqref{eq:OptStep5} has $T_f\times N_{in}^j$ integer decision variables constrained to be non-negative, and has $T_f\times N_j$ inequality constraints on traffic density.

\item Given $^jU_{t:t+T_f-1}^{t,\ast}$ as the optimal solution of \eqref{eq:OptStep5} obtained at time $t$, solve the following optimization problem to determine the configuration of traffic lights at intersection $I_j$ over the prediction horizon $[t,t+T_f]$ (i.e., $^j\lambda_{t:t+T_f-1}^{t,\ast}$):
\begin{subequations}\label{eq:OptStep6}
\begin{align}
\min\limits_{\lambda_j}\; \sum\limits_{k=0}^{T_f-1}\left\Vert ^j\tilde{x}_{df}(k|t)\right\Vert_{\Gamma_n^j}^2,
\end{align}
subject to
\begin{align}
^j\tilde{x}(k|t)\subseteq~^j\tilde{\mathcal{X}},&~k=1,\cdots,T_f,\\
\lambda_j(t+k)\in\Lambda_j,&~k=0,\cdots,T_f-1,
\end{align}
\end{subequations}
where  $^j\tilde{x}(k|t)$ can be computed via \eqref{eq:xtilde} with initial condition $x(t|t-1)$, and  $^j\tilde{\mathcal{X}}\subset\mathbb{R}_{\geq0}^{N_j}$ is a polyhedron containing the origin. The optimization problem \eqref{eq:OptStep6} has $T_f$ decision variables constrained to belong to the set $\Lambda_j$ (see Subsection \ref{sec:ActionSpace}), and has $T_f\times N_j$ inequality constraints on traffic density.

\end{enumerate}

Note that the above-mentioned scheme is receding horizon control-based; that is the Control Unit\#$j,~j\in\mathcal{M}$ computes the optimal inlet traffic flows and configuration of the traffic lights at intersection $I_j$ over the prediction horizon by solving the associated optimization problems at every time instant $t$, but only implements the next inlet traffic flows and configuration of traffic lights, and then solves the associated optimization problems again at the next time instant, repeatedly.

\begin{remark}
The optimization problem \eqref{eq:OptStep5} can be solved by mixed-integer tools, and the optimization problem \eqref{eq:OptStep6} can be solved by using the brute-force search algorithms.
\end{remark}

\begin{remark}
In decentralized scheme, Control Unit\#$j,~j\in\mathcal{M}$ estimates the traffic density at time $t$ across the network by assuming $d(t-1)=0$. Thus, in general, $x(t|t-1)\neq x(t)$. Also, Control Unit\#$j$ determines the optimal decisions over the prediction horizon based upon the optimal decisions of other control units at time $t-1$. As a result, the decentralized scheme is expected to provide a sub-optimal solution. This will be shown in Section \ref{sec:simulation}.
\end{remark}

\subsection{Emergency Traffic Mode}
Consider the assumptions mentioned in Subsection \ref{sec:EmergencyCondition} regarding the arriving, leaving, and recovery times. The control objective in an emergency traffic mode is to shorten the traveling time of the emergency vehicle, without causing congestion in other lanes. Given $T_a^t$, $T_s^t$, and $T_r^t$ by the aggregator, the Control Unit\#$j,~j\in\mathcal{M}$ executes the following steps to determine the inlet traffic flows and configuration of the traffic lights at intersection $I_j$ in an emergency traffic mode. Note that the following scheme is receding horizon control-based; that is the Control Unit\#$j,~j\in\mathcal{M}$ computes the optimal inlet traffic flows and configuration of the traffic lights at intersection $I_j$ over the prediction horizon by solving the associated optimization problems at every time instant $t$, but only implements the next inlet traffic flows and configuration of traffic lights, and then solves the associated optimization problems again at the next time instant, repeatedly.

\begin{enumerate}
\item Compute $x(t|t-1)$ based on the shared information by the aggregator, and according to \eqref{eq:system3} with $d(t-1)=0$.
\item Update traffic density at local lanes (i.e., $^jx(t)$), and replace corresponding elements in $x(t|t-1)$ with updated values.
\item Compute $\{^i\lambda_{t:t+T_f-2}^{t-1,\ast},\lambda_i(t+T_f-1)\}$ for all $i\in\mathcal{M}$, where $^i\lambda_{t:t+T_f-2}^{t-1,\ast}$ is the optimal solution of Control Unit\#$i$ obtained at time $t-1$ and $\lambda_i(t+T_f-1)$ is selected randomly from the action space $\Lambda_i$. 
\item Compute $\{^iU_{t:t+T_f-2}^{t-1,\ast},^iU_{nom}(t+T_f-1)\}$ for all $i\in\mathcal{M}$ and $i\neq j$, where $^iU_{t:t+T_f-2}^{t-1,\ast}$ is the optimal solution of Control Unit\#$i$ obtained at time $t-1$. 
\item Solve the following optimization problem to determine the inlet traffic flows at intersection $I_j$ over the prediction horizon $[t,t+T_f]$ (i.e., $^jU_{t:t+T_f-1}^{t,\ast}$):
\begin{subequations}\label{eq:OptStep7}
\begin{align}
\min\limits_{^jU}\; \sum\limits_{k=0}^{T_f-1}\Big(&\left\Vert ^j\hat{x}_{df}(k|t)\right\Vert_{^j\Gamma_e}^2\nonumber\\
&+\left\Vert ^jU(t+k)-~^jU_{nom}(t+k)\right\Vert_{^j\Theta}^2\Big),
\end{align}
subject to
\begin{align}
^j\hat{x}(k|t)\subseteq~^j\hat{\mathcal{X}}^+,&~k=1,\cdots,T_a^t+T_s^t+T_r^t,\\
^j\hat{x}(k|t)\subseteq~^j\hat{\mathcal{X}},&~k=T_a^t+T_s^t+T_r^t+1,\cdots,T_f,\\
^jU(t+k)\in\mathbb{Z}_{\geq0}^{N_{in}^j},&~k=0,\cdots,T_f-1,
\end{align}
\end{subequations}
where $^j\hat{\mathcal{X}}^+\supset~^j\hat{\mathcal{X}}$ is the extended set (see Remark \ref{remark:Extension}), and $^j\Gamma_e=~^j\Gamma_e^\top\geq0$ ($\in\mathbb{R}^{N_j}$) is the weighting matrix (see Remark \ref{remark:Gammae}). Similar to \eqref{eq:OptStep5}, the optimization problem \eqref{eq:OptStep7} has $T_f\times N_{in}^j$ integer decision variables constrained to be non-negative, and has $T_f\times N_j$ inequality constraints on traffic density.

\item Given $^jU_{t:t+T_f-1}^{t,\ast}$ as the optimal solution of \eqref{eq:OptStep7} obtained at time $t$, solve the following optimization problem to determine the configuration of traffic lights at intersection $I_j$ over the prediction horizon $[t,t+T_f]$ (i.e., $^j\lambda_{t:t+T_f-1}^{t,\ast}$):
\begin{subequations}\label{eq:OptStep8}
\begin{align}
\min\limits_{\lambda_j}\; \sum\limits_{k=0}^{T_f-1}\left\Vert ^j\tilde{x}_{df}(k|t)\right\Vert_{\Gamma_e^j}^2,
\end{align}
subject to
\begin{align}
^j\tilde{x}(k|t)\subseteq~^j\tilde{\mathcal{X}}^+,&~k=1,\cdots,T_a^t+T_s^t+T_r^t,\\
^j\tilde{x}(k|t)\subseteq~^j\tilde{\mathcal{X}},&~k=T_a^t+T_s^t+T_r^t+1,\cdots,T_f,\\
\lambda_j(t+k)\in\Lambda_j,&~k=0,\cdots,T_f-1,
\end{align}
\end{subequations}
where $^j\tilde{\mathcal{X}}^+\supset~^j\tilde{\mathcal{X}}$ is the extended set (see Remark \ref{remark:Extension}). Similar to \eqref{eq:OptStep6}, the optimization problem \eqref{eq:OptStep8} has $T_f$ decision variables constrained to belong to the set $\Lambda_j$ (see Subsection \ref{sec:ActionSpace}), and has $T_f\times N_j$ inequality constraints on traffic density.
\end{enumerate}

\begin{remark}
The optimization problem \eqref{eq:OptStep7} can be solved by mixed-integer tools, and the optimization problem \eqref{eq:OptStep8} can be solved by using the brute-force search algorithms.
\end{remark}

\begin{remark}
In decentralized scheme the emergency path $I_e^\ast$ is determined by the emergency vehicle, and is shared with control units through the aggregator.
\end{remark}

\begin{remark}
In this paper we assume that each control unit in the decentralized scheme controls the inlet traffic flows and configuration of traffic lights at one intersection. However, the decentralized scheme is applicable to the case where a network is divided into some sub-networks, and there exist a control unit in each sub-network controlling the entire sub-network.
\end{remark}

\section{Simulation Results}\label{sec:simulation}
\begin{figure}
\centering
\includegraphics[width=7cm]{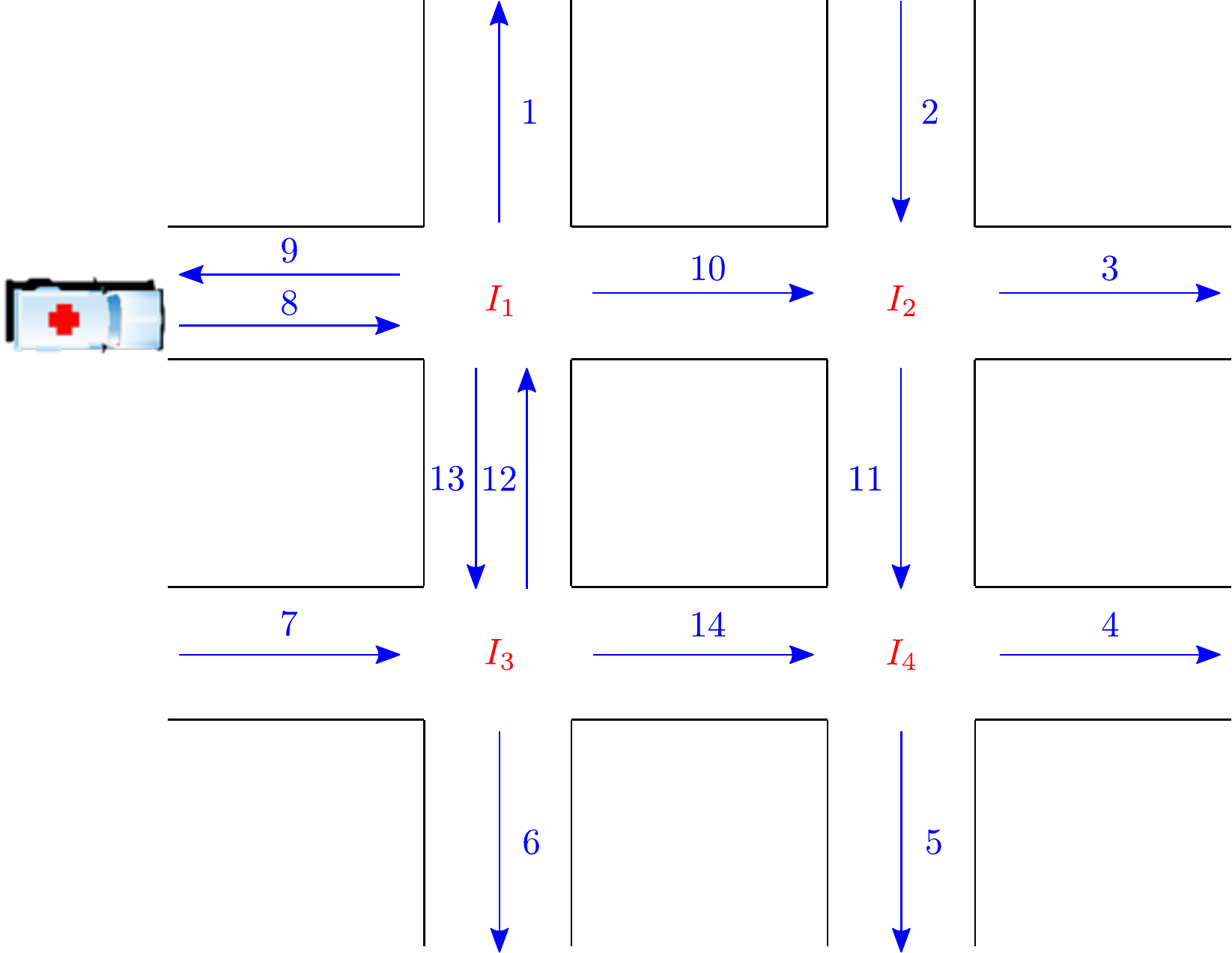}
\caption{Considered traffic network with 14 lanes and 4 intersections. An emergency vehicle enters through lane 8 and leaves through lane 5.}
\label{fig:Problem}
\end{figure}

Consider the traffic network shown in Fig. \ref{fig:Problem}. This network contains 14 unidirectional lanes identified by the set $\mathcal{N}=\{1,\cdots,14\}$, and 4 intersections identified by the set $\mathcal{M}=\{1,\cdots,4\}$. Also, $\mathcal{N}_{in}=\{2,7,8\}$. The edge set is $\mathcal{E}=\{(2,3),(2,11),(7,12),(7,14),(7,6),(8,1),(8,10),(8,13),\allowbreak(10,3),(10,11),(11,4),(11,5),(12,1),(12,9),(12,10),\allowbreak(13,6),(13,14),(14,4),(14,5)\}$.

\begin{figure}
\centering
\includegraphics[width=7cm]{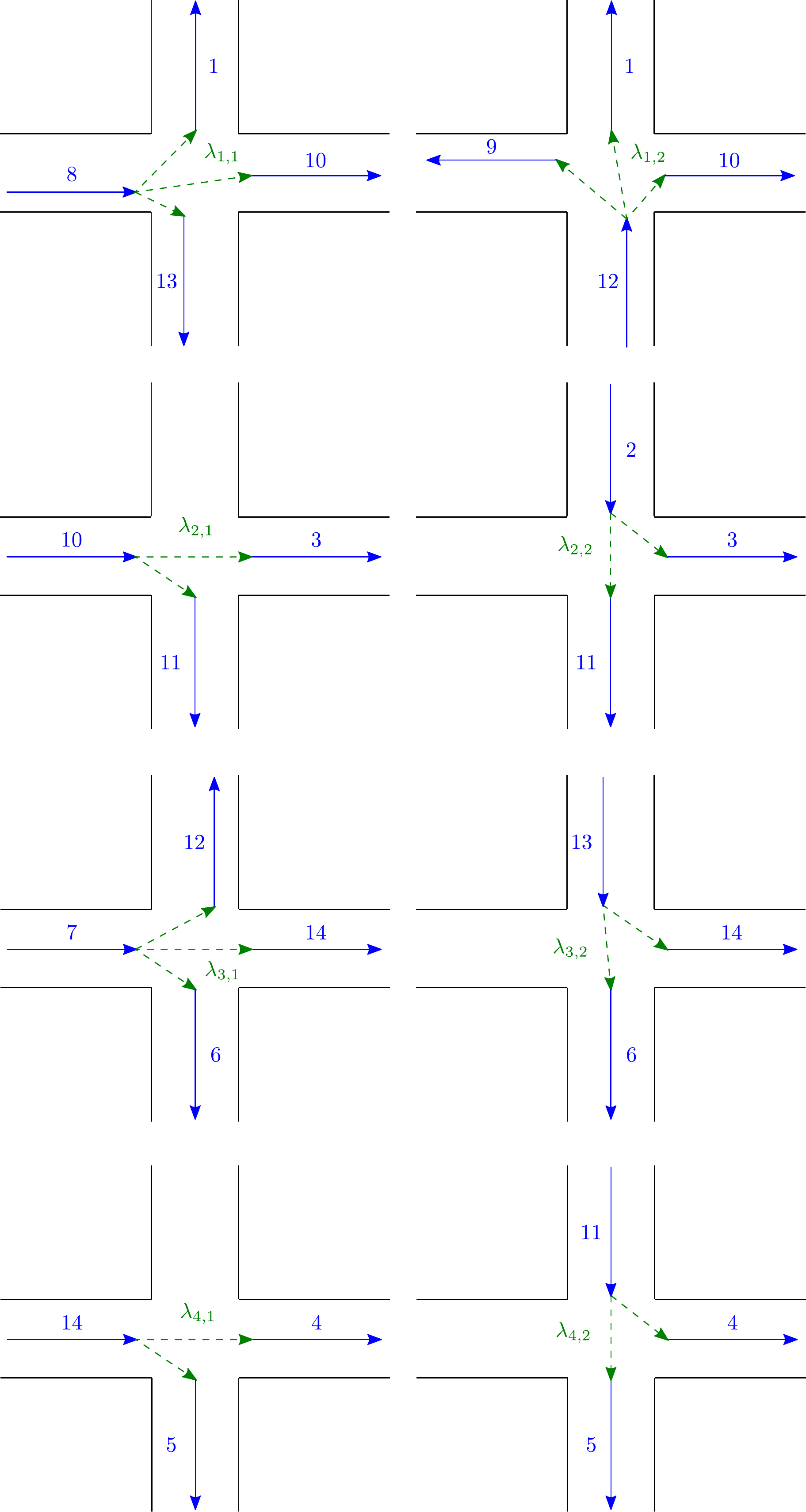}
\caption{Possible configurations of traffic lights at each intersection of the considered traffic network. The action space is $\Lambda=\Lambda_1\times\Lambda_2\times\Lambda_3\times\Lambda_4$, where $\Lambda_1=\{\lambda_{1,1},\lambda_{1,2}\}$, $\Lambda_2=\{\lambda_{2,1},\lambda_{2,2}\}$, $\Lambda_3=\{\lambda_{3,1},\lambda_{3,2}\}$, and $\Lambda_4=\{\lambda_{4,1},\lambda_{4,2}\}$.}
\label{fig:ActionSpace}
\end{figure}

Fig. \ref{fig:ActionSpace} shows possible configurations of traffic lights at each intersection of the traffic network shown in Fig. \ref{fig:Problem}. As seen in this figure, $\mu_1=\mu_2=\mu_3=\mu_4=2$, and the possible configurations at each intersection are: i) Intersection $I_1$: $\lambda_{1,1}$ corresponds to a `green' light at the end of lane 8, and a `red' light at the end of lane 12; $\lambda_{1,2}$ corresponds to a `red' light at the end of lane 8, and a `green' light at the end of lane 12; ii) Intersection $I_2$: $\lambda_{2,1}$ corresponds to a `green' light at the end of lane 10, and a `red' light at the end of lane 2; $\lambda_{2,2}$ corresponds to a `red' light at the end of lane 10, and a `green'  light at the end of lane 2; iii) Intersection $I_3$: $\lambda_{3,1}$ corresponds to a `green' light at the end of lane 7, and a `red' light at the end of lane 13; $\lambda_{3,2}$ corresponds to a `red' light at the end of lane 7, and a `green' light at the end of lane 13; and iv) Intersection $I_4$: $\lambda_{4,1}$ corresponds to a `green' light at the end of lane 14, and a `red' light at the end of lane 11; $\lambda_{4,2}$ corresponds to a `red' light at the end of lane 14, and a `green' light at the end of lane 11.

The boundary inflow vector of the traffic network shown in Fig. \ref{fig:Problem} is $U(t)=[u_2(t)~u_7(t)~u_8(t)]^\top\in\mathbb{Z}_{\geq0}^3$. We assume that $\Delta T=30$ seconds; this sampling period is appropriate to address macroscopic characteristics of traffic flow \cite{Rastgoftar2021,Erp2018,Wong2021}, as the traffic cycle ranges from one minute to three minutes in real-world systems \cite{NACTO}. For intersection $I_1$ and for the action $\lambda=(\lambda_{1,1},\lambda_2,\lambda_3,\lambda_4)$, we have $p_8(\lambda)\in[0,1]$, $p_{12}(\lambda)=0$, $q_{8,1}(\lambda),q_{8,10}(\lambda),q_{8,13}(\lambda)\in[0,1]$, and $q_{12,1}(\lambda)=q_{12,9}(\lambda)=q_{12,10}(\lambda)=0$. For intersection $I_1$ and for the action $\lambda=(\lambda_{1,2},\lambda_2,\lambda_3,\lambda_4)$, we have $p_8(\lambda)=0$, $p_{12}(\lambda)\in[0,1]$, $q_{8,1}(\lambda),q_{8,10}(\lambda),q_{8,13}(\lambda)=0$, and $q_{12,1}(\lambda)=q_{12,9}(\lambda)=q_{12,10}(\lambda)\in[0,1]$. For intersection $I_2$ and for the action $\lambda=(\lambda_1,\lambda_{2,1},\lambda_3,\lambda_4)$, we have $p_{10}(\lambda)\in[0,1]$, $p_{2}(\lambda)=0$, $q_{10,3}(\lambda),q_{10,11}(\lambda)\in[0,1]$, and $q_{2,3}(\lambda)=q_{2,11}(\lambda)=0$. For intersection $I_2$ and for the action $\lambda=(\lambda_1,\lambda_{2,2},\lambda_3,\lambda_4)$, we have $p_{10}(\lambda)=0$, $p_{2}(\lambda)\in[0,1]$, $q_{10,3}(\lambda),q_{10,11}(\lambda)=0$, and $q_{2,3}(\lambda)=q_{2,11}(\lambda)\in[0,1]$. For intersection $I_3$ and for the action $\lambda=(\lambda_1,\lambda_2,\lambda_{3,1},\lambda_4)$, we have $p_{13}(\lambda)=0$, $p_{7}(\lambda)\in[0,1]$, $q_{13,6}(\lambda),q_{13,14}(\lambda)=0$, and $q_{7,12}(\lambda)=q_{7,14}(\lambda)=q_{7,6}(\lambda)\in[0,1]$. For intersection $I_3$ and for the action $\lambda=(\lambda_1,\lambda_2,\lambda_{3,2},\lambda_4)$, we have $p_{13}(\lambda)\in[0,1]$, $p_{7}(\lambda)=0$, $q_{13,6}(\lambda),q_{13,14}(\lambda)\in[0,1]$, and $q_{7,12}(\lambda)=q_{7,14}(\lambda)=q_{7,6}(\lambda)=0$. For intersection $I_4$ and for the action $\lambda=(\lambda_1,\lambda_2,\lambda_3,\lambda_{4,1})$, we have $p_{14}(\lambda)\in[0,1]$, $p_{11}(\lambda)=0$, $q_{14,4}(\lambda),q_{14,5}(\lambda)\in[0,1]$, and $q_{11,4}(\lambda)=q_{11,5}(\lambda)=0$. For intersection $I_4$ and for the action $\lambda=(\lambda_1,\lambda_2,\lambda_3,\lambda_{4,2})$, we have $p_{14}(\lambda)=0$, $p_{11}(\lambda)\in[0,1]$, $q_{14,4}(\lambda),q_{14,5}(\lambda)=0$, and $q_{11,4}(\lambda)=q_{11,5}(\lambda)\in[0,1]$.

For implementing the decentralized scheme, we assume $^1x(t)=[x_1(t)~x_8(t)~x_9(t)~x_{12}(t)]^\top\in\mathbb{Z}_{\geq0}^4$, $^2x(t)=[x_2(t)~x_3(t)~x_{10}(t)]^\top\in\mathbb{Z}_{\geq0}^3$, $^3x(t)=[x_6(t)~x_7(t)~x_{13}(t)]^\top\in\mathbb{Z}_{\geq0}^3$, and $^4x(t)=[x_4(t)~x_5(t)~x_{11}(t)~x_{14}(t)]^\top\in\mathbb{Z}_{\geq0}^4$. That is Control Unit\#1 controls lanes 1, 8, 9, and 12; Control Unit\#2 controls lanes 2, 3, and 10; Control Unit\#3 controls lanes 6, 7, and 13; and Control Unit\#4 controls lanes 4, 5, 11, and 14. Also, $^1U(t)=u_8(t)$, $^2U(t)=u_2(t)$, and $^3U(t)=u_7(t)$. Thus, $N_{in}^1=1$, $N_{in}^2=1$, $N_{in}^3=1$, and $N_{in}^4=0$.

\begin{figure}[!t]
\centering
\includegraphics[width=8.5cm]{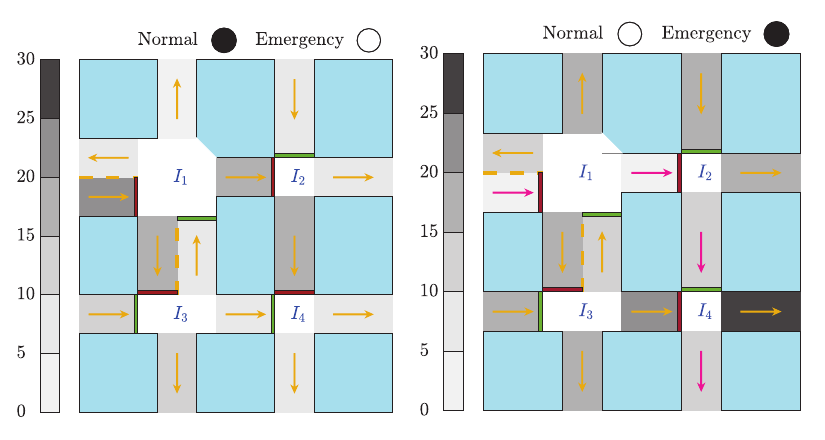}
\caption{A screenshot of the generated simulator shown in the accompanied video (\url{https://youtu.be/FmEYCxmD-Oc}). The black circle at the top shows the traffic mode in the network, which is either normal or emergency. The color of each lane indicates the traffic density, which can be interpreted according to the bar at the left. Yellow arrows show traffic direction at each lane, and pink arrows show selected emergency path.}
\label{fig:Simulator}
\end{figure}

The simulations are run on an Intel(R) Core(TM) i7-7500U CPU 2.70 GHz with 16.00 GB of RAM. In order to have a visual demonstration of the considered traffic network, a simulator is generated (see Fig. \ref{fig:Simulator}). A video of operation of the simulator is available at the URL: \url{https://youtu.be/FmEYCxmD-Oc}. For comparison purposes, we also simulate the centralized scheme presented in \cite{Rastgoftar2021} and a typical/existing/usual/baseline traffic system (i.e., the system with periodic schedule for traffic lights). TABLE \ref{tab:ComputationTime} compares the mean Computation Time (CT) of the proposed schemes per time step with that of the scheme presented in \cite{Rastgoftar2021}, where the value for the scheme of \cite{Rastgoftar2021} is used as the basis for normalization. As can be seen from this table, the computation time of the proposed centralized scheme is $\sim1.5$ times less than that of the scheme of \cite{Rastgoftar2021}. The computation time of the proposed decentralized scheme is $\sim1000$ times less than that of the scheme of \cite{Rastgoftar2021}, and is $\sim800$ times less than that of the proposed centralized scheme.

\begin{table}[!t]
\centering
\caption{Comparing the mean computation time of the proposed schemes with that of scheme of \cite{Rastgoftar2019}.}
\begin{tabular}{c|c|c|c}
& Centralized & Decentralized & Scheme of \cite{Rastgoftar2021} \\
\hline\hline
Mean CT (Norm.) & $0.734$ & $1.03\times10^{-3}$ & $1$
\end{tabular}
\label{tab:ComputationTime}
\end{table}

\subsection{Normal Traffic Mode}
Let $\hat{\mathcal{X}}=\tilde{\mathcal{X}}=\{x|x_i\leq20,~i\in\{1,\cdots,14\}\}$, and $d_i(t),~\forall i$ be selected uniformly from $\{-2,-1,0,1,2\}$. The initial condition is $x(0)=[15,16,15,12,12,17,18,10,10,14,12,10,16,10]^\top$, and the nominal boundary inflow is $U_{nom}(t)=[6,6,8]^\top$. Also, $\Theta=50I_{N_{in}}$ and $\gamma_i^n=1,~\forall i$.

\begin{figure}[!t]
\centering
\includegraphics[width=8.5cm]{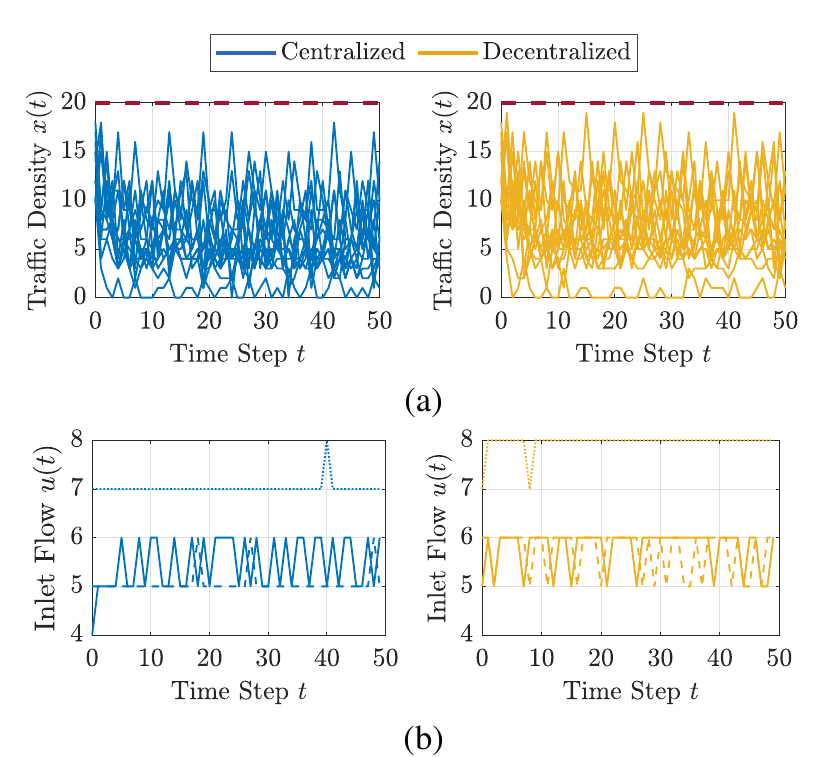}
\caption{Simulation results for the centralized and decentralized traffic control schemes in a normal traffic mode. Figure (a): Traffic density in all lanes. Figure (b): Inlet traffic flows.}
\label{fig:NormalCentralized}
\end{figure}

Simulation results are shown in Fig. \ref{fig:NormalCentralized}. TABLE \ref{tab:CentralizedNormal} compares the achieved Steady-State Density (SSD) with the considered schemes, where the value for the typical/existing/usual/baseline traffic system is used as the basis for normalization. Note that the reports are based on results of 1000 runs. According to TABLE \ref{tab:CentralizedNormal}, all methods perform better than the typical/existing/usual/baseline traffic system. The proposed centralized scheme provides the best response. The proposed decentralized scheme outperforms the scheme of \cite{Rastgoftar2021}, while as expected, it yields a larger SSD compared to the proposed centralized scheme. More precisely, degradation in the mean SSD by the decentralized scheme in comparison with the centralized scheme in a normal traffic mode is 11.42\% which is small and acceptable in real-life traffic scenarios. Thus, the cost of using the decentralized scheme instead of the centralized scheme in a normal traffic mode is very small.

\subsection{Emergency Traffic Mode}
Suppose that at time $t=10$, the aggregator receives a notification that an emergency vehicle will enter the network through lane 8 in two time steps, and should leave the network in two time steps through lane 5. Also, suppose that we have one time step to recover the traffic condition. We have $\mathcal{P}=\{I_e^1,I_e^2\}$, where $I_e^1=\{8,13,14,5\}$ and $I_e^2=\{8,10,11,5\}$. Let $\hat{\mathcal{X}}^+=\tilde{\mathcal{X}}^+=\{x|x_i\leq25\}$.

\begin{figure}[!t]
\centering
\includegraphics[width=8.5cm]{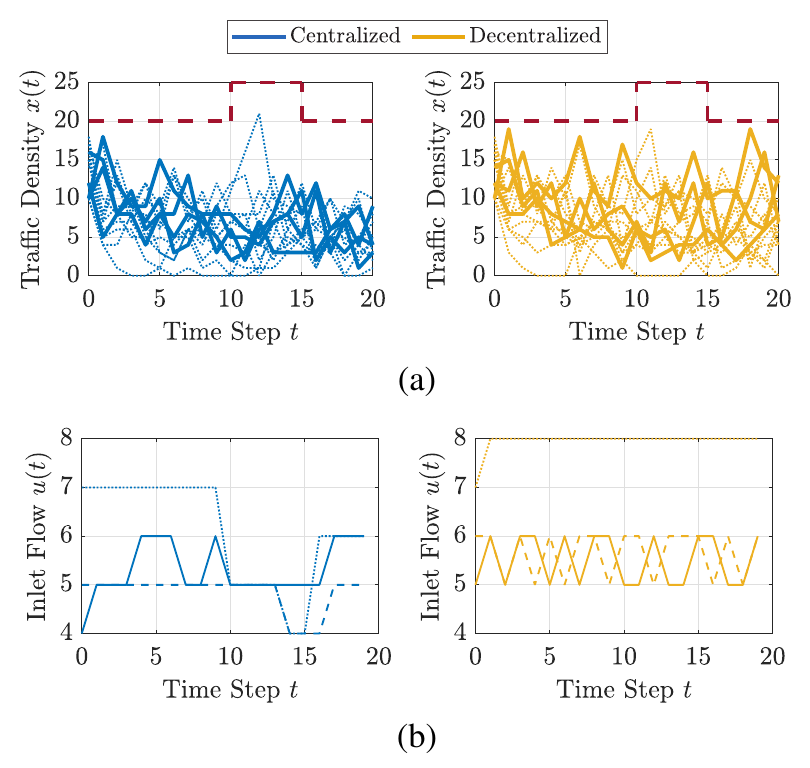}
\caption{Simulation results for the centralized and decentralized traffic control schemes in an emergency traffic mode. Figure (a): Traffic density in all lanes. Figure (b): Inlet traffic flows.}
\label{fig:EmergencyDistributed}
\end{figure}

Simulation results are shown in Fig. \ref{fig:EmergencyDistributed}, and the results from comparison analysis are reported in TABLE \ref{tab:CentralizedEmergency}, that are computed based on results of 1000 runs. Note that the values for the typical/existing/usual/baseline  traffic system are used as nominal values for normalization. As seen in TABLE \ref{tab:CentralizedEmergency}, both schemes proposed in this paper perform better than the typical/existing/usual/baseline traffic system in an emergency traffic mode. In particular, the centralized and decentralized schemes reduce the mean SSD by 23.73\% and 14.58\%, respectively. As expected, the decentralized scheme yields a larger SSD compared to the centralized scheme. More precisely, degradation in mean SSD by the decentralized scheme in comparison with the centralized scheme is 11.98\%. TABLE \ref{tab:CentralizedEmergency} also reports that the centralized and decentralized schemes reduce the mean Density in Emergency Path (DEP) by 47.97\% and 30.42\%, respectively. It is noteworthy that the degradation in the mean DEP by the decentralized scheme in comparison with the centralized scheme is 33.73\%.

\begin{table}[!t]
\centering
\caption{Comparing the mean steady-state density with the considered schemes in a normal traffic mode. The value for the typical/existing/usual/baseline traffic system is used as the basis for normalization.}
\begin{tabular}{c|c|c|c}
& Centralized & Decentralized & Scheme of \cite{Rastgoftar2021} \\
\hline\hline
Mean SSD (Norm.) & $0.7251$ & $0.8079$ & $0.9035$
\end{tabular}
\label{tab:CentralizedNormal}
\end{table}

\begin{table}[!t]
\centering
\caption{Comparing the mean steady-state density and density in emergency path with the proposed centralized and decentralized schemes in an emergency traffic mode. The values for the typical/existing/usual/baseline traffic system are used as bases for normalization.}
\begin{tabular}{c|c|c}
& Centralized & Decentralized  \\
\hline\hline
Mean SSD (Norm.) & 0.7627 & 0.8542 \\
Mean DEP (Norm.) & 0.5203 & 0.6958
\end{tabular}
\label{tab:CentralizedEmergency}
\end{table}

\subsection{Sensitivity Analysis\textemdash Impact of Look-Ahead Horizon $T_f$}
\begin{figure}
\centering
\includegraphics[width=8.5cm]{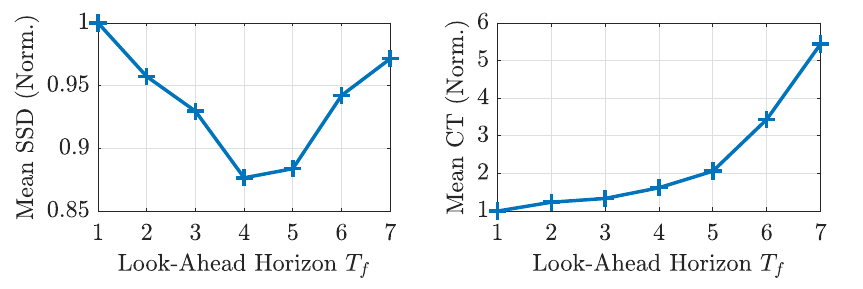}
\caption{Impact of the look-ahead horizon $T_f$ on the performance and computation time of the proposed decentralized scheme, where the values for $T_f=1$ are used as bases for normalization. As $T_f$ increases the performance improves at the cost of increased computation time. For large $T_f$ the performance degrades due to poor prediction accuracy.}
\label{fig:sensitivity}
\end{figure}

Fig. \ref{fig:sensitivity} shows how the prediction window size impacts the performance and computation time of the developed centralized scheme, where the values for $T_f=1$ are used as nominal values for normalization. From Fig. \ref{fig:sensitivity}-left, we see that as the look-ahead horizon increases, the performance of the decentralized scheme improves as it takes into account more information of future conditions. However, as we look further into the future, the performance is degraded since prediction accuracy reduces. Fig. \ref{fig:sensitivity}-right shows that as the look-ahead horizon increases, the computation time of the proposed decentralized scheme increases concomitantly with the size and complexity of the associated optimization problems. In the simulation studies, we selected $T_f=4$, as it yields the best performance with an affordable computing time. Note that a similar behavior is observed for the centralized scheme; that is $T_f=4$ provides the best performance for the centralized scheme.

\section{Conclusion}\label{sec:conclusion}
This paper proposed an emergency vehicle-centered traffic control framework to alleviate traffic congestion in a network of interconnected signaled lanes. The aim of this paper is to integrate CTM with MPC, to ensure that emergency vehicles traverse multiple intersections efficiently and timely. Two schemes were developed in this paper: i) centralized; and ii) decentralized. It was shown that the centralized scheme provides the optimal solution, even though its computation time may be large for large networks. To cope with this problem, a decentralized scheme was developed, where an aggregator acts as the hub of the network. It was shown that the computation time of the decentralized scheme is very small, which makes it a good candidate for large networks, even though it provides a sub-optimal solution. Extensive simulation studies were carried out to validate and evaluate the performance of the proposed schemes.

Future work will aim at extending the developed schemes to deal with cases where two (or more) emergency vehicles traverse a network. This extension is not trivial, and requires addressing many technical and methodological challenges. Also, future work should investigate robustness and tolerance of the decentralized scheme to uncertainty in communication delay and communication failures.

\bibliographystyle{IEEEtran}
\bibliography{ref}{}

\begin{IEEEbiography}[{\includegraphics[width=1in,height=1.25in,clip,keepaspectratio]{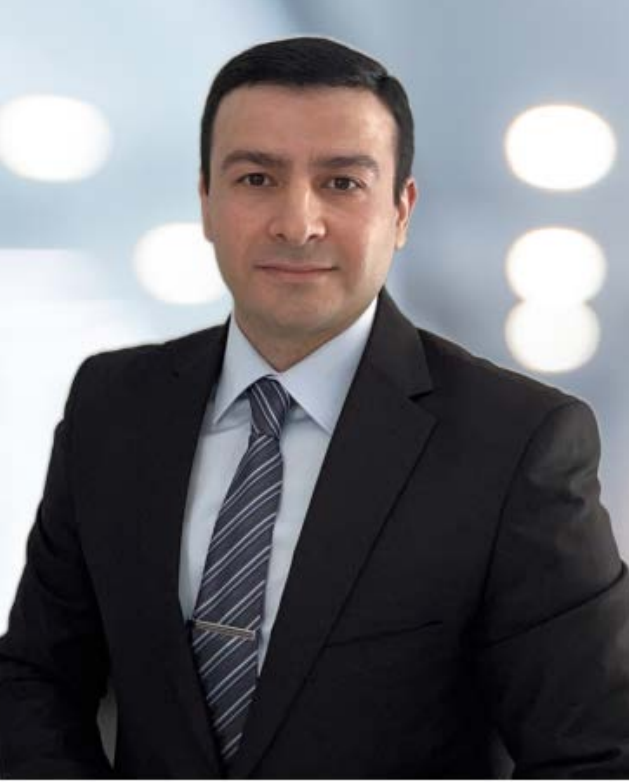}}] {Mehdi Hosseinzadeh} received his Ph.D. degree in Electrical Engineering-Control from the University of Tehran, Iran, in 2016. From 2017 to 2019, he was a postdoctoral researcher at Universit\'{e} Libre de Bruxelles, Brussels, Belgium. In 2018, he was a visiting researcher at University of British Columbia, Canada. He is currently a postdoctoral research associate at Washington University in St. Louis, MO, USA. His research interests include nonlinear and adaptive control, constrained control, and safe and robust control of autonomous systems. \vspace{-4cm}
\end{IEEEbiography}

\begin{IEEEbiography}[{\includegraphics[width=1in,height=1.25in,clip,keepaspectratio]{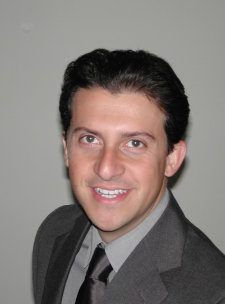}}] {Bruno  Sinopoli} received his Ph.D. in Electrical Engineering from the University of  California at Berkeley, in 2005. After a postdoctoral position at Stanford University, he was the faculty at Carnegie Mellon University from 2007 to 2019, where he was full  professor  in  the  Department  of  Electrical  and Computer  Engineering  with  courtesy  appointments in  Mechanical  Engineering  and  in  the  Robotics  Institute  and  co-director  of  the  Smart  Infrastructure Institute. In  2019  he  joined Washington University in Saint Louis, where he is the chair of the Electrical and Systems Engineering department. He was awarded the 2006 Eli Jury  Award  for  outstanding  research  achievement  in  the  areas  of  systems, communications,  control  and  signal  processing  at  U.C.  Berkeley,  the  2010 George Tallman Ladd Research Award from Carnegie Mellon University and the NSF Career award in 2010. His research interests include the modeling,analysis  and  design  of  Secure  by  Design  Cyber-Physical  Systems  with applications  to  Energy  Systems,  Interdependent  Infrastructures  and  Internet of Things. \vspace{-4cm}
\end{IEEEbiography}

\begin{IEEEbiography}[{\includegraphics[width=1in,height=1.25in,clip,keepaspectratio]{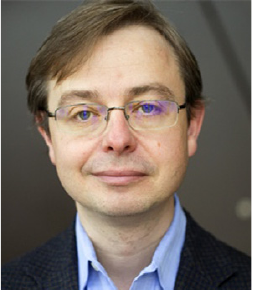}}] {Ilya Kolmanovsky}  is a professor in the department of aerospace engineering at the University of Michigan,with research interests in control theory for systems with state and control constraints, and in control applications to aerospace and automotive systems. He received his Ph.D. degree in aerospace engineering from the University of Michigan in 1995. Prof. Kolmanovsky is a Fellow of IEEE and is named as an inventor on 104 United States patents. \vspace{-4cm}
\end{IEEEbiography}

\begin{IEEEbiography}[{\includegraphics[width=1in,height=1.25in,clip,keepaspectratio]{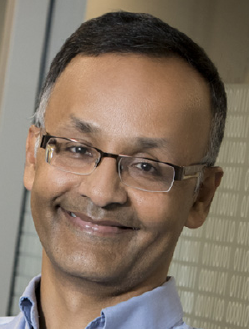}}] {Sanjoy Baruah} joined Washington University in St. Louis in September 2017. He was previously at the University of North Carolina at Chapel Hill (1999\textendash 2017) and the University of Vermont (1993\textendash 1999). His research interests and activities are in real-time and safety-critical system design, scheduling theory, resource allocation and sharing in distributed computing environments, and algorithm design and analysis.
\end{IEEEbiography}

\end{document}